\documentclass[12pt]{article}

\usepackage{fullpage}
\usepackage{amsmath,amssymb,amsthm}
\usepackage{tikz}
\usetikzlibrary{arrows.meta,calc}
\usepackage{booktabs}
\usepackage{microtype}
\usepackage[colorlinks=true,linkcolor=blue,citecolor=blue,urlcolor=blue]{hyperref}

\newtheorem{theorem}{Theorem}[section]
\newtheorem{proposition}[theorem]{Proposition}
\newtheorem{lemma}[theorem]{Lemma}
\newtheorem{corollary}[theorem]{Corollary}
\theoremstyle{definition}
\newtheorem{example}[theorem]{Example}
\theoremstyle{remark}
\newtheorem{remark}[theorem]{Remark}

\newcommand{\partitions}{\mathcal{P}}
\newcommand{\E}{\mathbb{E}}
\newcommand{\Var}{\operatorname{Var}}
\newcommand{\seqnum}[1]{\href{https://oeis.org/#1}{#1}}

\begin{document}

\begingroup
\renewcommand{\thefootnote}{\fnsymbol{footnote}}
\begin{center}
\vskip 1cm
{\LARGE\bfseries Binomially Weighted Partition Statistics:\\[4pt]
Generating Functions, Asymptotics, and Dominance}
\vskip 0.85cm
\large
Grant Kopitzke, Alex Phillippi\footnotemark[2],
Derek Scheiderer\footnotemark[2] \\
Department of Mathematics\\
University of Wisconsin--Stevens Point\\
Stevens Point, WI 54481, USA\\
\href{mailto:gkopitzk@uwsp.edu}{\texttt{gkopitzk@uwsp.edu}}
\end{center}
\footnotetext[2]{Undergraduate student}
\endgroup
\vskip .2 in

\begin{abstract}
For an integer partition $\lambda=(\lambda_1,\lambda_2,\ldots)$, we consider
$M_m(\lambda)=\sum_{i\geq1}\binom{i+m}{m+1}\lambda_i$, where $m\geq0$.
We prove conjugation and Frobenius-coordinate identities for this family and
obtain its ordinary and bivariate generating functions. Differentiating the
bivariate product gives exact formulas for the mean and variance under
uniform sampling of partitions of $n$. For fixed $m$, we show that
$\mathbb E[M_m]\sim\zeta(m+3)(\sqrt{6n}/\pi)^{m+3}$ and determine an
explicit positive constant multiplying $(\sqrt{6n}/\pi)^{2m+5}$ in the
variance asymptotic. The family encompasses the row sum underlying the
discrete Gini statistic and, after conjugation, a linear combination giving
the Wiener index of starlike trees; the associated normalized statistics
have contrasting limits under uniform random sampling. We also prove
that every $M_m$ is strictly decreasing along the dominance order.
Consequently its level sets are antichains, providing coefficient-based
lower bounds for the width of the dominance lattice.
\end{abstract}

\section{Introduction}\label{sec:introduction}

An \emph{integer partition} of a positive integer $n$ is a weakly decreasing sequence $\lambda=(\lambda_1,\ldots,\lambda_\ell)$ of positive integers with $\sum_{i=1}^{\ell}\lambda_i=n$. We write $\lambda\vdash n$, and call $\ell=\ell(\lambda)$ its \emph{length}. The empty partition $\varnothing$ is the unique partition of $0$. For example, the partitions of $4$ are
\[
(4),\quad(3,1),\quad(2,2),\quad(2,1,1),\quad(1,1,1,1).
\]
We use $p(n)$ for the number of partitions of $n$, with $p(0)=1$, and
$\partitions$ for the set of partitions of all nonnegative integers.
Standard background on partitions and Young diagrams can be found in
\cite{Andrews:1,Fulton:1}.

Two apparently different questions motivate our investigation. The
\emph{discrete Gini statistic} of an integer partition is a measure
inspired by inequality and Lorenz curves, while the \emph{Wiener index}
of a starlike tree measures the total distance between pairs of
vertices. Both can be expressed using the same elementary family of
binomially weighted sums, with partition conjugation accounting for the
different ways that the two objects are naturally encoded. This common
structure makes it possible to study both examples, as well as higher
members of the family, through one set of combinatorial and analytic tools.

The mechanism behind this connection is a simple observation about Young
diagrams: a statistic obtained by summing a weight over all cells can be
evaluated either row by row or column by column. Choosing binomial weights
is especially effective because sums of consecutive binomial coefficients
collapse by the hockey-stick identity. The resulting column formulas are
amenable to partition generating functions, while the row formulas reveal
a monotonicity phenomenon in the dominance order. Thus the same construction
links symmetric-function identities, random partitions, and antichains.

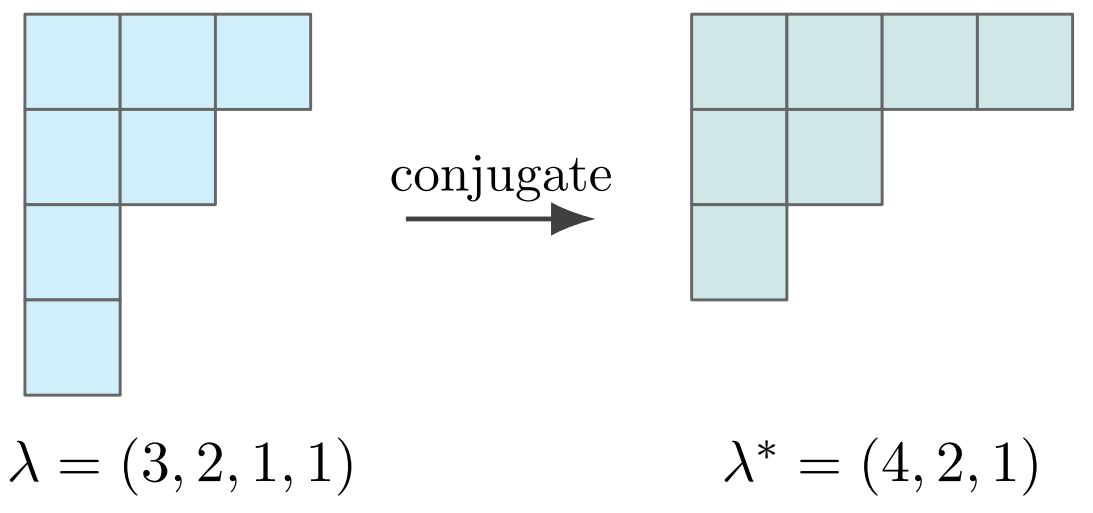
\begin{figure}[htbp]
\centering
\begin{tikzpicture}[scale=0.63, line join=round]
  \begin{scope}[shift={(0,0)}]
    \foreach \r/\len in {1/3,2/2,3/1,4/1}{
      \foreach \c in {1,...,\len}{
        \filldraw[fill=cyan!16,draw=black!60,line width=0.55pt]
          ({\c-1},{1-\r}) rectangle (\c,{2-\r});
      }
    }
    \node[font=\small] at (1.65,-3.75) {$\lambda=(3,2,1,1)$};
  \end{scope}
  \draw[-{Latex[length=3mm]},line width=0.9pt,black!75]
    (4.0,-1.15)--(6.0,-1.15);
  \node[font=\footnotesize] at (5.0,-0.72) {conjugate};
  \begin{scope}[shift={(7,0)}]
    \foreach \r/\len in {1/4,2/2,3/1}{
      \foreach \c in {1,...,\len}{
        \filldraw[fill=teal!19,draw=black!60,line width=0.55pt]
          ({\c-1},{1-\r}) rectangle (\c,{2-\r});
      }
    }
    \node[font=\small] at (2.0,-3.75) {$\lambda^*=(4,2,1)$};
  \end{scope}
\end{tikzpicture}
\caption{Young diagrams for a partition and its conjugate. Each column
of $Y(\lambda)$ becomes a row of $Y(\lambda^*)$.}
\label{fig:conjugation}
\end{figure}

The \emph{Young diagram} $Y(\lambda)$ consists of left-justified cells,
with $\lambda_i$ cells in row $i$. The \emph{conjugate partition}
$\lambda^*$ is obtained by reflecting $Y(\lambda)$ across its main
diagonal. Equivalently,
\[
 \lambda_j^*=\bigl|\{i:\lambda_i\geq j\}\bigr|,
 \qquad (3,2,1,1)^*=(4,2,1).
\]
The diagrams in Figure~\ref{fig:conjugation} illustrate this correspondence.
In particular, conjugation exchanges row lengths with column heights while
preserving the total number of cells, so it is a bijection on the partitions
of each fixed integer. We will use this bijection repeatedly: some of our
formulas are naturally expressed using row indices, while the generating
functions become simpler when expressed using multiplicities of parts of
the conjugate partition.

Partitions are padded with trailing zeros whenever an expression requires
more parts than their lengths.

An important example of a statistic on partitions is the unnormalized
discrete Gini statistic of \cite{Kopitzke:1}, motivated by the Gini
index and Lorenz curve in the study of inequality
\cite{Origins:1,Lorenz:1,Farris:1}. In the convention used here,
\begin{equation}\label{gini}
 g(\lambda)=\binom{n+1}{2}-\sum_{i\geq1}i\lambda_i.
\end{equation}
The elementary conjugation identity proved in
Proposition~\ref{prop:basecase} shows that
\begin{equation}\label{gini=e2}
 g(\lambda)=\binom{n+1}{2}
   -\sum_{j=1}^{\lambda_1}\binom{\lambda_j^*+1}{2}
   =e_2(\lambda^*),
\end{equation}
where $e_2(\nu)=\sum_{i<j}\nu_i\nu_j$ is the second elementary
symmetric polynomial evaluated at the parts of $\nu$.
The statistic in \eqref{gini} should not be confused with its normalized
version, obtained by dividing by $\binom n2$ for $n\geq2$.
The equality with $e_2(\lambda^*)$ is useful because it connects the
inequality interpretation of $g$ with the elementary symmetric functions
of the column heights. For example, $(n)$ has conjugate $(1^n)$ and hence
$g((n))=\binom n2$, whereas $(1^n)$ has conjugate $(n)$ and $g((1^n))=0$.
Thus the normalization takes values in $[0,1]$ and its orientation is
determined by the convention in \eqref{gini}.

A related expression appears in the Wiener index of a starlike tree.
For a connected graph $T$, its \emph{Wiener index} is
$W(T)=\sum_{\{u,v\}\subseteq V(T)}d_T(u,v)$, the sum of graph
distances over unordered pairs of distinct vertices. For example,
$W(K_{1,3})=3\cdot1+3\cdot2=9$.
Given a partition $\lambda\vdash n$, construct a tree $S_\lambda$
by attaching paths of lengths $\lambda_1,\ldots,\lambda_\ell$ to a
common vertex. When $\ell\geq3$, the resulting tree has a unique
vertex of degree at least $3$, and its unrooted isomorphism type is
uniquely determined by the parts of $\lambda$. For $\ell\leq2$ the
construction gives a path, and different partitions may give the same
unrooted tree. The formula discussed in \cite{WienerIndex:1} is
\begin{equation}\label{wiener=e3}
 W(S_\lambda)=\binom{n+2}{3}-e_3(\lambda),
 \qquad e_3(\lambda)=\sum_{i<j<k}\lambda_i\lambda_j\lambda_k.
\end{equation}
To interpret this formula, recall that an edge of a tree separates its
vertices into sets of sizes $j$ and $n+1-j$, and contributes
$j(n+1-j)$ to the Wiener index: precisely those pairs with one vertex
on each side use that edge in their unique connecting path. On an arm of
length $\lambda_i$, these sizes are $j$ and $n+1-j$ for
$1\leq j\leq\lambda_i$. Consequently,
\[
 W(S_\lambda)
 =\sum_i\sum_{j=1}^{\lambda_i}j(n+1-j).
\]
The elementary power-sum identities reduce this expression to
\eqref{wiener=e3}. We give a related calculation in
Proposition~\ref{prop:wienerWeightedRelation} below. The formulas \eqref{gini=e2} and \eqref{wiener=e3} therefore place
an inequality statistic and a graph-distance statistic side by side:
although the former involves $e_2(\lambda^*)$ and the latter involves
$e_3(\lambda)$, both ultimately arise from weighted sums over rows or
columns of Young diagrams. The change between $\lambda$ and
$\lambda^*$ must be retained in any pointwise comparison.

These examples motivate a family of binomially weighted partition sums.
For $m\in\mathbb Z_{\geq0}$, put
\begin{equation}\label{eq:defineMm}
 M_m(\lambda)=\sum_{i=1}^{\ell(\lambda)}
               \binom{i+m}{m+1}\lambda_i,
 \qquad M_m(\varnothing)=0.
\end{equation}
For $m=0$, the weight is $i$, and for $m=1$ the weights on
successive rows are $\binom{i+1}{2}=1,3,6,10,\ldots$. More generally,
$\binom{i+m}{m+1}$ is a polynomial of degree $m+1$ in the row index.
The direct connection to the motivating statistics is
\begin{equation}\label{eq:introGiniWienerLink}
 \begin{aligned}
  g(\lambda)&=\binom{n+1}{2}-M_0(\lambda),\\
  W(S_\lambda)&=(n+2)M_0(\lambda^*)-2M_1(\lambda^*).
 \end{aligned}
\end{equation}
The first equality is immediate from the definition of $g$; the
second is proved in Proposition~\ref{prop:wienerWeightedRelation}.
Equation~\eqref{eq:introGiniWienerLink} explains why the family is useful
beyond the binomial identities themselves. The Gini statistic is an
affine transform of its first member, whereas the Wiener index combines
the first two members evaluated on the conjugate partition. Both
examples therefore inherit generating-function information and
probabilistic consequences from results for $M_m$. The conjugation in
the Wiener expression is essential for individual partitions, even
though it disappears when averaging over all partitions of a fixed
integer.

As $m$ varies, the weights give a hierarchy of partition statistics
that extends these two examples rather than merely reproducing them.
Our purpose is to examine the hierarchy through column summation,
principal hooks, generating functions, uniform random partitions,
and dominance order.

\medskip
\noindent\textbf{Contributions and relation to earlier work.}
The $m=0$ statistic and its connection with the Gini index were
considered in \cite{Kopitzke:1,Kopitzke:Dissertation}, while the Wiener
index of starlike trees and its symmetric-function expression have
been studied in \cite{WienerIndex:1}. Standard ingredients such as
Young-diagram double counting, Euler products, and conjugation of
dominance are not claimed as new in isolation. The contributions of
the present paper are the following systematic results for the whole
binomially weighted family.

First, Proposition~\ref{prop:generalcase} and
Proposition~\ref{prop:frobenius} give conjugation and principal-hook
descriptions that connect row weights to part-size statistics. The
product formula in Proposition~\ref{prop:twoVariableGF} extends the
$bivariate$ Gini generating function from the case $m=0$ to all
$m\geq0$, and Theorem~\ref{theorem:moments} derives exact first and
second moments from this product. Second,
Theorem~\ref{theorem:asymptoticMoments} establishes leading
expectation and variance asymptotics, with an explicit strictly
positive variance constant. In particular, the calculation of the
variance requires a second-order saddle-point correction because
the dominant second-moment terms cancel. The applications in
Corollaries~\ref{cor:normalizedGiniAsymptotic} and
\ref{cor:normalizedWienerAsymptotic} show a contrasting phenomenon:
for the normalizations used here, the Gini index tends to $1$ in
probability whereas the Wiener index tends to $0$. Finally,
Proposition~\ref{prop:dominancelattice} establishes strict
monotonicity under dominance and yields antichains indexed by values
of $M_m$, giving computable lower bounds for lattice width.

The paper is organized accordingly. Section~\ref{sec:identities}
develops Young-diagram and Frobenius-coordinate identities and proves
the Wiener relation in \eqref{eq:introGiniWienerLink}.
Section~\ref{sec:generatingFunctions} treats the one-variable and
bivariate generating functions, exact moments, and asymptotics.
Section~\ref{sec:dominanceAntichains} applies the same statistics to
dominance and antichains, and Section~\ref{sec:conclusion} summarizes
further directions.

\section{Weighted Partition Identities}\label{sec:identities}
\subsection{Conjugate-partition identities}

\begin{proposition}\label{prop:basecase}
For every partition $\lambda\vdash n$,
\begin{equation}\label{eq:basecase}
 \sum_{i=1}^{\ell(\lambda)}i\lambda_i
 =\sum_{j=1}^{\lambda_1}\binom{\lambda_j^*+1}{2}.
\end{equation}
\end{proposition}

The preceding identity is the $m=0$ case of the following more general
row--column summation identity. We include a combinatorial proof to
make the weighting convention explicit.
The key point is that the weight is attached to a \emph{row index}, not
to the length of the row. Conjugation changes these row-indexed weights
into functions of column heights, and these are the quantities that later
appear in the generating functions.

\begin{proposition}\label{prop:generalcase}
For every partition $\lambda\vdash n$ and every $m\in\mathbb Z_{\geq0}$,
\begin{equation}\label{formula:generalcase}
 M_m(\lambda)=
 \sum_{j=1}^{\lambda_1}\binom{\lambda_j^*+m+1}{m+2}.
\end{equation}
\end{proposition}

\begin{proof}
Assign the weight $w_i=\binom{i+m}{m+1}$ to every cell in row
$i$ of $Y(\lambda)$. There are $\lambda_i$ cells in that row, so
summing the weights by rows gives
\[
 \sum_{i=1}^{\ell(\lambda)}\lambda_iw_i=M_m(\lambda).
\]
Instead, fix a column $j$. Its height is $\lambda_j^*$, so it meets
exactly rows $1,2,\ldots,\lambda_j^*$. The same cells contribute
\[
 \sum_{i=1}^{\lambda_j^*}w_i
 =\sum_{i=1}^{\lambda_j^*}\binom{i+m}{m+1}
 =\binom{\lambda_j^*+m+1}{m+2}.
\]
For the last equality, apply Pascal's identity
$\binom{h+1}{r+1}-\binom{h}{r+1}=\binom{h}{r}$ and telescope over
$h=m+1,\ldots,\lambda_j^*+m$. Because the row sum and column sum
count exactly the same weighted cells, summing the displayed expression
over all columns proves \eqref{formula:generalcase}.
\end{proof}

\begin{proof}[Proof of Proposition~\ref{prop:basecase}]
Take $m=0$ in Proposition~\ref{prop:generalcase}.
\end{proof}

More generally, for any function $f:\mathbb Z_{>0}\to\mathbb C$,
counting row weights in the same Young diagram by columns gives
\begin{equation}\label{eq:arbitraryRowWeights}
 \sum_{i=1}^{\ell(\lambda)}f(i)\lambda_i
 =\sum_{j=1}^{\lambda_1}\sum_{i=1}^{\lambda_j^*}f(i).
\end{equation}
Proposition~\ref{prop:generalcase} follows by choosing
$f(i)=\binom{i+m}{m+1}$.

\begin{example}\label{ex:conjugateWeights}
Let $\lambda=(3,2,1,1)$, so $\lambda^*=(4,2,1)$.
The row-index formulation for $m=0$ gives
\[
 M_0(\lambda)=1\cdot3+2\cdot2+3\cdot1+4\cdot1=14.
\]
The column-height formulation gives the same answer:
\[
 \binom{4+1}{2}+\binom{2+1}{2}+\binom{1+1}{2}
 =10+3+1=14.
\]
For $m=1$, the same comparison becomes
\[
 M_1(\lambda)=3+3\cdot2+6+10=25
 =\binom63+\binom43+\binom33.
\]
This example emphasizes that
$\sum_i\binom{\lambda_i+m+1}{m+2}=M_m(\lambda^*)$,
not $M_m(\lambda)$ in general.
\end{example}

\subsection{Frobenius notation}\label{sec:frobenius}

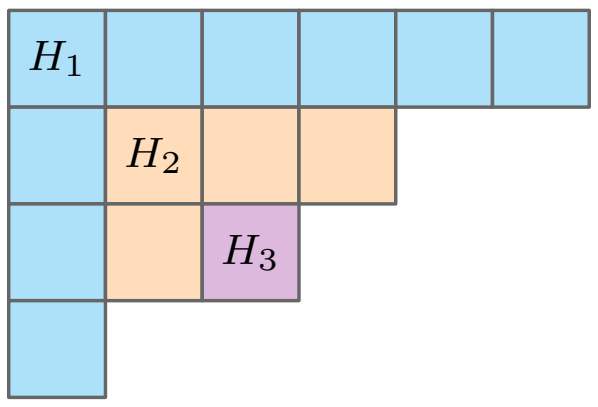
\begin{figure}[htbp]
\centering
\begin{tikzpicture}[scale=0.64, line join=round]
  \foreach \r/\len in {1/6,2/4,3/3,4/1}{
    \foreach \c in {1,...,\len}{
      \pgfmathtruncatemacro{\hookindex}{min(\r,\c)}
      \ifcase\hookindex\relax
      \or \def\hookcolor{cyan!29}
      \or \def\hookcolor{orange!27}
      \or \def\hookcolor{violet!27}
      \fi
      \filldraw[fill=\hookcolor,draw=black!60,line width=0.65pt]
        ({\c-1},{1-\r}) rectangle (\c,{2-\r});
      \ifnum\r=\c
         \node[font=\small,inner sep=0pt] at ({\c-0.5},{1.5-\r})
           {$\scriptstyle H_{\hookindex}$};
      \fi
    }
  }
  \node[font=\small,anchor=west] at (6.5,0.35)
    {$\lambda=(6,4,3,1)$};
  \node[font=\small,anchor=west] at (6.5,-0.50)
    {$d=3$};
  \node[font=\small,anchor=west] at (6.5,-1.35)
    {$\alpha=(5,2,0)$};
  \node[font=\small,anchor=west] at (6.5,-2.20)
    {$\beta=(3,1,0)$};
\end{tikzpicture}
\caption{The three principal hooks of the partition $(6,4,3,1)$,
shown in different colors. Hook $H_i$ consists of cell $(i,i)$, its arm
to the right, and its leg below. Their disjoint union is the diagram.}
\label{fig:principalHooks}
\end{figure}

The \emph{Durfee size} $d$ of $\lambda$ is the number of cells on
its main diagonal, equivalently $d=\max\{i:\lambda_i\geq i\}$.
Its \emph{Frobenius coordinates} are
\[
 \lambda=(\alpha_1,\ldots,\alpha_d\mid\beta_1,\ldots,\beta_d),
 \qquad \alpha_i=\lambda_i-i,\quad\beta_i=\lambda_i^*-i.
\]
Thus $\alpha_i$ and $\beta_i$ count, respectively, the cells strictly
to the right of and strictly below the diagonal cell $(i,i)$.
This notation is standard in partition theory; see, for example,
\cite{Andrews:1,Fulton:1}. Figure~\ref{fig:principalHooks} illustrates
the decomposition into principal hooks used in the next proof.
The principal hook based at $(i,i)$ contains its diagonal cell, the
$\alpha_i$ cells extending to the right, and the $\beta_i$ cells
extending downward. Every cell belongs to exactly one such hook: a cell
$(r,c)$ belongs to the hook based at $(\min\{r,c\},\min\{r,c\})$.
This partitions the diagram into disjoint pieces, making it possible to
compute the same row-weighted statistic without explicitly listing all
rows or columns.

\begin{proposition}\label{prop:frobenius}
Let $\lambda\vdash n$ have Frobenius coordinates
$(\alpha_1,\ldots,\alpha_d\mid\beta_1,\ldots,\beta_d)$.
For $m\in\mathbb Z_{\geq0}$,
\begin{equation}\label{eq:frobenius}
 M_m(\lambda)
 =\sum_{i=1}^d\left[
   \binom{i+m}{m+1}\alpha_i
   +\sum_{j=i}^{i+\beta_i}\binom{j+m}{m+1}
  \right].
\end{equation}
\end{proposition}

\begin{proof}
As in the proof of Proposition~\ref{prop:generalcase}, give every
cell in row $j$ of $Y(\lambda)$ the weight $\binom{j+m}{m+1}$.
For each diagonal cell $(i,i)$, let $H_i$ be its \emph{principal
hook}: the diagonal cell together with all cells strictly to its right
and strictly below it. The $d$ principal hooks are pairwise disjoint
and together contain every cell of $Y(\lambda)$.

The arm of $H_i$ contains $\alpha_i$ cells in row $i$, each of
weight $\binom{i+m}{m+1}$; its total contribution is therefore
$\alpha_i\binom{i+m}{m+1}$. The diagonal cell is in row $i$, and
the leg has one cell in each row $i+1,\ldots,i+\beta_i$. Their
combined contribution is
$\sum_{j=i}^{i+\beta_i}\binom{j+m}{m+1}$.
Since every cell belongs to one and only one principal hook, the sum
of these contributions over $i=1,\ldots,d$ is precisely $M_m(\lambda)$.
This is \eqref{eq:frobenius}.
\end{proof}

\begin{corollary}\label{cor:frobeniusClosed}
Under the hypotheses of Proposition~\ref{prop:frobenius},
\begin{equation}\label{eq:frobeniusClosed}
 M_m(\lambda)
 =\sum_{i=1}^d\left[
 \binom{i+m}{m+1}\alpha_i
 +\binom{i+\beta_i+m+1}{m+2}
 -\binom{i+m}{m+2}
 \right].
\end{equation}
\end{corollary}
\begin{proof}
For each $i$, the hockey-stick identity gives
\[
 \sum_{j=i}^{i+\beta_i}\binom{j+m}{m+1}
 =\binom{i+\beta_i+m+1}{m+2}-\binom{i+m}{m+2}.
\]
Substitute into \eqref{eq:frobenius}.
\end{proof}

\begin{proposition}\label{prop:wienerWeightedRelation}
For every partition $\lambda\vdash n$, the Wiener index of its
associated starlike tree satisfies
\begin{equation}\label{eq:wienerWeightedRelation}
 W(S_\lambda)=(n+2)M_0(\lambda^*)-2M_1(\lambda^*).
\end{equation}
\end{proposition}
\begin{proof}
Write $S_j=\sum_i\lambda_i^j$, so $S_1=n$. By the conjugation
identity \eqref{formula:generalcase},
\[
 M_0(\lambda^*)=\sum_i\binom{\lambda_i+1}{2}
 =\frac{S_2+n}{2},
 \qquad
 M_1(\lambda^*)=\sum_i\binom{\lambda_i+2}{3}
 =\frac{S_3+3S_2+2n}{6}.
\]
On the other hand, expanding the elementary symmetric polynomial
gives Newton's identity
\[
 6e_3(\lambda)=S_1^3-3S_1S_2+2S_3
 =n^3-3nS_2+2S_3.
\]
Using \eqref{wiener=e3}, we consequently obtain
\begin{align*}
 W(S_\lambda)
 &=\frac{n(n+1)(n+2)}6
   -\frac{n^3-3nS_2+2S_3}{6}\\
 &=\frac{3nS_2+3n^2+2n-2S_3}{6}\\
 &=(n+2)\frac{S_2+n}{2}
      -2\frac{S_3+3S_2+2n}{6}\\
 &=(n+2)M_0(\lambda^*)-2M_1(\lambda^*).
\end{align*}
\end{proof}

The conjugation in \eqref{eq:wienerWeightedRelation} is essential.
For example, when $n\geq3$ and $\lambda=(n)$, its associated tree is a path, but
$M_0(\lambda)=M_1(\lambda)=n$; substituting $\lambda$ instead of
$\lambda^*$ would incorrectly give $n^2$. Nevertheless,
conjugation preserves the uniform distribution on partitions of $n$,
so the distinction disappears after taking expectations. This will
be useful in Section~\ref{sec:moments}.

\section{Generating Functions}\label{sec:generatingFunctions}

We work with formal power series in $x$; every coefficient in the
infinite sums and products below depends on only finitely many terms.
An ordinary generating function records information according to the
size $|\lambda|$ of a partition, whereas a second marking variable
will later retain the actual values of the statistic $M_m$.
A product formula is particularly natural here because a partition
is uniquely specified by the multiplicity of each possible part size.
Recall Euler's generating function for the partition function
\cite{Andrews:1,Wilf:1}:
\begin{equation}\label{eq:EulerProduct}
 P(x)=\sum_{n\geq0}p(n)x^n=\prod_{k\geq1}\frac1{1-x^k}.
\end{equation}
Indeed, expanding $(1-x^k)^{-1}=\sum_{r\geq0}x^{kr}$ records the
multiplicity of each part size $k$.

\subsection{One-variable generating functions}\label{sec:onevariableGF}

Let $g_m(n)=\sum_{\lambda\vdash n}M_m(\lambda)$, and set
$G_m(x)=\sum_{n\geq1}g_m(n)x^n$.
Thus $[x^n]G_m(x)$ is the \emph{sum} of the statistic over all
partitions of $n$, not the number of partitions having a prescribed
statistic value. Dividing this coefficient by $p(n)$ will later give
its expectation under the uniform distribution.

\begin{theorem}\label{theorem:generalGF}
For every $m\in\mathbb Z_{\geq0}$,
\begin{equation}\label{eq:generalGF}
 G_m(x)=P(x)\sum_{k\geq1}\frac{x^k}{(1-x^k)^{m+3}}.
\end{equation}
\end{theorem}

For a partition $\lambda$, let $m_k(\lambda)$ be the number of
parts equal to $k$. Then $|\lambda|=\sum_{k\geq1}k\,m_k(\lambda)$.
We first record a standard multiplicity-generating-function identity.
The ordinary Euler factor $(1-x^k)^{-1}$ records how often the part
$k$ occurs. To sum those multiplicities over all partitions, we mark
the corresponding factor with an auxiliary variable and differentiate
before returning the marker to $1$.

\begin{lemma}\label{lemma:multiplicityGF}
For every $k\geq1$,
\begin{equation}\label{eq:multiplicityGF}
 F_k(x):=\sum_{\lambda\in\partitions}m_k(\lambda)x^{|\lambda|}
 =\frac{x^k}{1-x^k}P(x).
\end{equation}
\end{lemma}
\begin{proof}
Introduce a marker $u$ for each occurrence of a part of size $k$.
Choosing $r$ copies of $k$ contributes $u^rx^{kr}$, so the
corresponding geometric series is $1/(1-ux^k)$. The other part sizes
are unrestricted, and consequently
\[
 \sum_{\lambda\in\partitions}u^{m_k(\lambda)}x^{|\lambda|}
 =\frac1{1-ux^k}\prod_{j\ne k}\frac1{1-x^j}.
\]
Differentiate formally with respect to $u$ and set $u=1$:
\[
 F_k(x)=\frac{x^k}{(1-x^k)^2}\prod_{j\ne k}\frac1{1-x^j}
       =\frac{x^k}{1-x^k}P(x).
\]
\end{proof}
The derivative works because $\left.\frac{d}{du}u^r\right|_{u=1}=r$:
it converts a marker for the multiplicity into that multiplicity itself.

The following Lambert-series identity permits us to pass between
weights attached to part sizes and a more compact form of the
generating function. It is a reindexing identity, not a new partition
bijection: both sides encode the same divisor sums.\begin{lemma}\label{lemma:lambert2}
For $m\in\mathbb Z_{\geq0}$,
\begin{equation}\label{eq:generalLambert}
 \sum_{k\geq1}\binom{k+m+1}{m+2}\frac{x^k}{1-x^k}
 =\sum_{k\geq1}\frac{x^k}{(1-x^k)^{m+3}}.
\end{equation}
\end{lemma}
\begin{proof}
Put $a_k=\binom{k+m+1}{m+2}$. The coefficient of $x^n$ in the
left-hand side is $\sum_{k\mid n}a_k$, since
$\frac{x^k}{1-x^k}=\sum_{j\geq1}x^{kj}$.
On the other hand, the negative-binomial expansion gives
\[
 \frac{x^k}{(1-x^k)^{m+3}}
 =\sum_{j\geq1}\binom{j+m+1}{m+2}x^{kj}
 =\sum_{j\geq1}a_jx^{kj}.
\]
Thus the coefficient of $x^n$ on the right is
$\sum_{k\mid n}a_{n/k}=\sum_{k\mid n}a_k$, because
$k\mapsto n/k$ permutes the divisors of $n$.
\end{proof}

\begin{corollary}\label{lemma:lambert1}
\begin{equation}\label{eq:baseLambert}
 \sum_{k\geq1}\binom{k+1}{2}\frac{x^k}{1-x^k}
 =\sum_{k\geq1}\frac{x^k}{(1-x^k)^3}.
\end{equation}
\end{corollary}
\begin{proof}
Take $m=0$ in Lemma~\ref{lemma:lambert2}.
\end{proof}

\begin{proof}[Proof of Theorem~\ref{theorem:generalGF}]
By Proposition~\ref{prop:generalcase}, grouping equal parts of the
conjugate gives
\begin{equation}\label{eq:generalMultiplicityWeight}
 M_m(\lambda)
 =\sum_{k\geq1}\binom{k+m+1}{m+2}m_k(\lambda^*).
\end{equation}
Conjugation preserves size and is a bijection of $\partitions$;
hence Lemma~\ref{lemma:multiplicityGF} yields
\begin{align*}
 G_m(x)
 &=\sum_{\lambda\in\partitions}M_m(\lambda)x^{|\lambda|}\\
 &=\sum_{k\geq1}\binom{k+m+1}{m+2}
       \sum_{\lambda\in\partitions}m_k(\lambda^*)x^{|\lambda|}\\
 &=P(x)\sum_{k\geq1}\binom{k+m+1}{m+2}\frac{x^k}{1-x^k}\\
   &=P(x)\sum_{k\geq1}\frac{x^k}{(1-x^k)^{m+3}},
\end{align*}
where the last step uses Lemma~\ref{lemma:lambert2}.
\end{proof}
The conjugation step is indispensable. We sum
$\binom{k+m+1}{m+2}$ over parts of size $k$ of $\lambda^*$,
not of $\lambda$ itself; only after reindexing the full sum over all
partitions can we replace $\lambda^*$ by a new partition variable.

For clarity, the $m=0$ case is recorded separately.
\begin{proposition}\label{prop:baseGF}
With $g(n)=\sum_{\lambda\vdash n}\sum_i i\lambda_i$ and
$G(x)=\sum_{n\geq1}g(n)x^n$,
\begin{equation}\label{eq:baseGF}
 G(x)=P(x)\sum_{k\geq1}\frac{x^k}{(1-x^k)^3}.
\end{equation}
\end{proposition}
\begin{proof}
Take $m=0$ in Theorem~\ref{theorem:generalGF}.
\end{proof}

\subsection{Two-variable generating functions}\label{sec:twovariableGF}

The one-variable series $G_m(x)$ records the total weight of
partitions of each $n$, but not the distribution of those weights.
Introducing a variable $q$ to mark $M_m(\lambda)$ retains the number
of partitions with each possible statistic value. Distinct partitions
can nevertheless have the same weight.
More precisely, $[q^rx^n]H_m(q,x)$ will count the partitions of
$n$ for which $M_m(\lambda)=r$. This additional information is needed
both for moment computations and for constructing antichains in the
dominance lattice. Since $M_m$ is a sum of weights over parts of the
conjugate partition, each occurrence of a part of size $k$ in that
conjugate contributes the factor $q^{\binom{k+m+1}{m+2}}x^k$.

\begin{proposition}\label{prop:twoVariableGF}
For $m\in\mathbb Z_{\geq0}$, let
\[
 H_m(q,x)=\sum_{n\geq0}\sum_{\lambda\vdash n}
                  q^{M_m(\lambda)}x^n.
\]
Then
\begin{equation}\label{eq:twoVariableGF}
 H_m(q,x)=\prod_{k\geq1}
   \frac1{1-q^{\binom{k+m+1}{m+2}}x^k}.
\end{equation}
\end{proposition}

The $m=0$ case appears in \cite{Kopitzke:1}. We include the general
proof to emphasize the role of conjugation.
\begin{proof}
Let $a_k=\binom{k+m+1}{m+2}$. Expanding the product gives
\begin{align*}
 \prod_{k\geq1}\frac1{1-q^{a_k}x^k}
 &=\prod_{k\geq1}\left(\sum_{r\geq0}q^{ra_k}x^{kr}\right)\\
 &=\sum_{\mu\in\partitions}
   q^{\sum_{k\geq1}a_km_k(\mu)}x^{|\mu|}.
\end{align*}
By \eqref{eq:generalMultiplicityWeight},
$\sum_k a_km_k(\mu)=M_m(\mu^*)$.
Conjugation is a size-preserving bijection, so this is
$\sum_{\lambda\in\partitions}q^{M_m(\lambda)}x^{|\lambda|}$.
The empty partition accounts for the constant term $1$.
\end{proof}

\begin{example}\label{ex:twoVariableSmall}
For $m=0$ and $n=4$, the five partitions listed in
Section~\ref{sec:introduction} have $M_0$ values
$4,5,6,7,10$, respectively. Thus the coefficient of $x^4$ is
\[
 [x^4]H_0(q,x)=q^4+q^5+q^6+q^7+q^{10}.
\]
Setting $q=1$ gives $p(4)=5$; applying $q\frac{\partial}{\partial q}$
and then setting $q=1$ gives $4+5+6+7+10=32$. The next
subsection develops this observation systematically.
\end{example}

\begin{corollary}\label{cor:recoverGeneralGF}
\begin{equation}\label{eq:recoverGeneralGF}
 H_m(1,x)=P(x),\qquad
 \left.\frac{\partial H_m(q,x)}{\partial q}\right|_{q=1}=G_m(x).
\end{equation}
\end{corollary}
\begin{proof}
The first identity follows by setting $q=1$. For the second,
formal differentiation gives
\[
 \left.\frac{\partial H_m(q,x)}{\partial q}\right|_{q=1}
 =\sum_{n\geq0}\sum_{\lambda\vdash n}M_m(\lambda)x^n
 =G_m(x).
\]
\end{proof}

\subsection{Moments under the uniform distribution}\label{sec:moments}

Fix $n\geq1$ and choose a partition $\Lambda_n$ uniformly from the
$p(n)$ partitions of $n$. The random variable
$X_{m,n}=M_m(\Lambda_n)$ takes the value $M_m(\lambda)$ with probability
$1/p(n)$ for each $\lambda\vdash n$. In particular,
\[
 \E[X_{m,n}]=\frac1{p(n)}\sum_{\lambda\vdash n}M_m(\lambda),
 \qquad
 \E[X_{m,n}^2]=\frac1{p(n)}\sum_{\lambda\vdash n}M_m(\lambda)^2.
\]
Our two-variable generating function contains all these values at
once: the exponent of $x$ records the size of the partition, and the
exponent of $q$ records the statistic. Differentiation with respect to
$q$ therefore converts the marked exponent into a multiplicative
weight. This avoids enumerating all partitions separately for each $n$.

For concision, define
\begin{equation}\label{eq:momentSeries}
 a_k=\binom{k+m+1}{m+2},\qquad
 L_m(x)=\sum_{k\geq1}a_k\frac{x^k}{1-x^k},\qquad
 Q_m(x)=\sum_{k\geq1}a_k^2\frac{x^k}{(1-x^k)^2}.
\end{equation}
By Lemma~\ref{lemma:lambert2},
$L_m(x)=\sum_{k\geq1}x^k/(1-x^k)^{m+3}$.
The square on $a_k$ in $Q_m$ will arise from taking two derivatives;
it should not be confused with the square of the series $L_m(x)$.

\begin{theorem}\label{theorem:moments}
For $n\geq1$ and $m\in\mathbb Z_{\geq0}$,
\begin{align}
 \E[X_{m,n}]
 &=\frac{[x^n]P(x)L_m(x)}{p(n)},
 \label{eq:expectation}\\
 \E[X_{m,n}^2]
 &=\frac{[x^n]P(x)\bigl(L_m(x)^2+Q_m(x)\bigr)}{p(n)},
 \label{eq:secondMoment}\\
 \Var(X_{m,n})
 &=\frac{[x^n]P(x)\bigl(L_m(x)^2+Q_m(x)\bigr)}{p(n)}
   -\left(\frac{[x^n]P(x)L_m(x)}{p(n)}\right)^2.
 \label{eq:variance}
\end{align}
\end{theorem}

\begin{proof}
Use the differential operator $D_q=q\,\partial/\partial q$. For a
nonnegative integer $h$, we have
\[
 D_q(q^h)=h q^h,\qquad D_q^2(q^h)=h^2q^h.
\]
Thus, on applying $D_q$ to
$H_m(q,x)=\sum_{n\geq0}\sum_{\lambda\vdash n}
q^{M_m(\lambda)}x^n$ and setting $q=1$, the coefficient of $x^n$
becomes $\sum_{\lambda\vdash n}M_m(\lambda)$. Applying $D_q$
twice instead gives $\sum_{\lambda\vdash n}M_m(\lambda)^2$.
It remains to calculate the two derivatives of the product
\eqref{eq:twoVariableGF}.

Because $H_m(q,x)$ has constant term $1$ as a series in $x$, its
formal logarithm is defined. The product formula gives
\[
 \log H_m(q,x)=-\sum_{k\geq1}\log(1-q^{a_k}x^k).
\]
Differentiating the summand $-\log(1-q^{a_k}x^k)$ yields
\[
 D_q\bigl[-\log(1-q^{a_k}x^k)\bigr]
 =\frac{a_kq^{a_k}x^k}{1-q^{a_k}x^k}.
\]
A second application of $D_q$, using the quotient rule, gives
\[
 D_q\left(\frac{a_kq^{a_k}x^k}{1-q^{a_k}x^k}\right)
 =\frac{a_k^2q^{a_k}x^k}{(1-q^{a_k}x^k)^2}.
\]
Consequently,
\begin{align*}
 \left.D_q\log H_m(q,x)\right|_{q=1}
 &=\sum_{k\geq1}a_k\frac{x^k}{1-x^k}=L_m(x),\\
 \left.D_q^2\log H_m(q,x)\right|_{q=1}
 &=\sum_{k\geq1}a_k^2\frac{x^k}{(1-x^k)^2}=Q_m(x).
\end{align*}

To recover derivatives of $H_m$ rather than its logarithm, use
$D_qH_m=H_mD_q\log H_m$. Applying the product rule once more,
\[
 D_q^2H_m
 =H_m\left((D_q\log H_m)^2+D_q^2\log H_m\right).
\]
Since $H_m(1,x)=P(x)$, it follows that
\begin{align*}
 \left.D_qH_m(q,x)\right|_{q=1}&=P(x)L_m(x),\\
 \left.D_q^2H_m(q,x)\right|_{q=1}
 &=P(x)\bigl(L_m(x)^2+Q_m(x)\bigr).
\end{align*}
Notice that the term $L_m(x)^2$ arises from the product rule;
omitting it would give an incorrect second moment. Extracting the
coefficient of $x^n$ and dividing by $p(n)$ proves
\eqref{eq:expectation} and \eqref{eq:secondMoment}. Finally,
\eqref{eq:variance} is the identity
$\Var(X)=\E[X^2]-\E[X]^2$.
\end{proof}

The coefficient formulas can be converted into finite sums involving
divisors. Indeed,
\[
 \frac{x^k}{1-x^k}=\sum_{j\geq1}x^{kj},
\]
so $x^r$ appears in the $k$th summand precisely when $k$ divides
$r$. For $r\geq1$, define
\begin{equation}\label{eq:bmDef}
 b_m(r)=\sum_{d\mid r}\binom{d+m+1}{m+2}.
\end{equation}
It follows that $L_m(x)=\sum_{r\geq1}b_m(r)x^r$. The product
with $P(x)=\sum_{j\geq0}p(j)x^j$ is then a Cauchy product: to
obtain $x^n$, the two exponents must add to $n$.

\begin{corollary}\label{cor:expectedConvolution}
For $n\geq1$,
\begin{equation}\label{eq:expectedConvolution}
 \E[X_{m,n}]=\frac1{p(n)}\sum_{r=1}^n p(n-r)b_m(r).
\end{equation}
In particular, if $\sigma_s(r)=\sum_{d\mid r}d^s$, then
\begin{equation}\label{eq:baseExpectation}
 \E[X_{0,n}]=\frac1{2p(n)}\sum_{r=1}^n
 p(n-r)\bigl(\sigma_1(r)+\sigma_2(r)\bigr).
\end{equation}
\end{corollary}
\begin{proof}
The Cauchy product gives
\[
 [x^n]P(x)L_m(x)
 =\sum_{r=1}^n p(n-r)b_m(r),
\]
proving the first claim by \eqref{eq:expectation}. If $m=0$,
then $\binom{d+1}{2}=(d^2+d)/2$, and hence
\[
 b_0(r)=\frac12\sum_{d\mid r}(d^2+d)
 =\frac12\bigl(\sigma_2(r)+\sigma_1(r)\bigr).
\]
Substitution into \eqref{eq:expectedConvolution} proves the
second identity.
\end{proof}

There is a parallel divisor-sum expansion for $Q_m$. The geometric
series identity
\[
 \frac{y}{(1-y)^2}=\sum_{j\geq1}jy^j
\]
shows that $x^k/(1-x^k)^2=\sum_{j\geq1}jx^{kj}$. Thus the
coefficient of $x^r$ is obtained by choosing a divisor $d=k$ of
$r$, at which point $j=r/d$. Put
\[
 c_m(r)=\sum_{d\mid r}\frac{r}{d}
     \binom{d+m+1}{m+2}^{\!2}.
\]
Then $Q_m(x)=\sum_{r\geq1}c_m(r)x^r$. Likewise,
$[x^r]L_m(x)^2=\sum_{u=1}^{r-1}b_m(u)b_m(r-u)$,
because both factors of $L_m$ have zero constant term. A second
application of the Cauchy product therefore yields
\begin{equation}\label{eq:secondMomentConvolution}
 \E[X_{m,n}^2]=\frac1{p(n)}\sum_{r=1}^n p(n-r)
 \left(c_m(r)+\sum_{u=1}^{r-1}b_m(u)b_m(r-u)\right).
\end{equation}
These expressions are exact finite formulas and are convenient for
numerical checks without multiplying the infinite products directly.

We next connect the moment identities with the Gini and Wiener
statistics from the introduction. For comparison with
\cite{Kopitzke:Dissertation}, define
$\widetilde g(\lambda)=g(\lambda)/\binom n2$ for
$\lambda\vdash n\geq2$. By \eqref{gini},
$g(\lambda)=\binom{n+1}{2}-M_0(\lambda)$, so
\[
 \E[\widetilde g(\Lambda_n)]
 =\frac{n+1}{n-1}-\frac{\E[X_{0,n}]}{\binom n2}.
\]
Inserting \eqref{eq:baseExpectation} gives
\begin{equation}\label{eq:expectedNormalizedGini}
 \E[\widetilde g(\Lambda_n)]
 =\frac{n+1}{n-1}
  -\frac{\sum_{r=1}^n p(n-r)
         (\sigma_1(r)+\sigma_2(r))}{n(n-1)p(n)},
\end{equation}
recovering the exact expectation formula in that work.

For the Wiener index, there is a choice of normalization. Dividing by
$\binom{n+2}{3}$ compares the Wiener index with its maximum, but
does not make its minimum zero. Here we use instead the \emph{min--max
normalization}. Among trees with $n$ edges, the star has minimum
Wiener index $n^2$ and the path has maximum Wiener index
$\binom{n+2}{3}$. Accordingly, for $n\geq3$ define
\[
 \widetilde W(\lambda)
 =\frac{W(S_\lambda)-n^2}{\binom{n+2}{3}-n^2}.
\]
The denominator is $\binom n3$, and \eqref{wiener=e3} simplifies
this to
\[
 \widetilde W(\lambda)=1-\frac{e_3(\lambda)}{\binom n3}.
\]
It takes values in $[0,1]$, with $0$ at $(1^n)$ (the star) and $1$
at $(n)$ (a path). The normalized Gini statistic has the same endpoint orientation
under our convention, although the two statistics behave differently
for a typical uniformly sampled partition.

For a pointwise identity we must use the conjugate partition, as in
Proposition~\ref{prop:wienerWeightedRelation}. Because
$\lambda\mapsto\lambda^*$ is a bijection on partitions of $n$,
uniform sampling gives
\[
 \E[W(S_{\Lambda_n})]
 =(n+2)\E[X_{0,n}]-2\E[X_{1,n}].
\]
We may therefore apply \eqref{eq:expectedConvolution} to both
terms. Since
\[
 b_0(r)=\frac{\sigma_1(r)+\sigma_2(r)}2,
 \qquad
 b_1(r)=\frac{\sigma_3(r)+3\sigma_2(r)+2\sigma_1(r)}6,
\]
we obtain
\begin{equation}\label{eq:expectedNormalizedWiener}
 \E[\widetilde W(\Lambda_n)]
 =\frac{\displaystyle\sum_{r=1}^n p(n-r)
 \left((3n+2)\sigma_1(r)+3n\sigma_2(r)-2\sigma_3(r)\right)
 -6n^2p(n)}{n(n-1)(n-2)p(n)}.
\end{equation}
The conjugation does not disappear at the level of individual
partitions; it disappears only after averaging over all partitions.
This distinction is important when applying the weighted-sum formulas
to graph invariants.

\begin{example}\label{ex:momentsN4}
For $n=4$, the values of $M_0$ on the five partitions, in the order
listed in Section~\ref{sec:introduction}, are $4,5,6,7,10$.
Their sum is $32$ and the sum of their squares is
$16+25+36+49+100=226$. Thus
\[
 \E[X_{0,4}]=\frac{32}{5},\qquad
 \E[X_{0,4}^2]=\frac{226}{5},\qquad
 \Var(X_{0,4})=\frac{226}{5}-\left(\frac{32}{5}\right)^2
 =\frac{106}{25}.
\]
For comparison, the five values of the Wiener index are
$20,20,20,18,16$; subtracting the minimum $16$ and dividing by
$\binom43=4$ gives normalized values $1,1,1,\tfrac12,0$.
Hence $\E[\widetilde W(\Lambda_4)]=7/10$.
\end{example}

\subsection{Asymptotic expectation and variance}
\label{sec:asymptotics}

The preceding formulas determine the moments exactly, but extracting
their coefficients becomes less practical as $n$ grows. We therefore
ask how the expectation and variance behave as $n\to\infty$ when the
weight parameter $m$ is fixed. The principal advantage of the
generating-function approach is that this question can be answered by
studying the singular behavior of $L_m(x)$ and $Q_m(x)$ near $x=1$,
where Euler's partition generating function is concentrated during
coefficient extraction.

Our analysis combines a saddle-point estimate for Euler's partition
product with Mellin expansions of the Lambert series
\cite{FlajoletSedgewick:1,Fristedt:1}. We write
$\zeta(z)=\sum_{j\geq1}j^{-z}$ for $\Re z>1$, with its standard
meromorphic continuation, and $\Gamma(z)$ for the gamma function.
The latter satisfies $\Gamma(k)=(k-1)!$ for positive integers $k$.
The notation $f(n)\sim g(n)$ means $f(n)/g(n)\to1$.

\begin{theorem}\label{theorem:asymptoticMoments}
Fix $m\in\mathbb Z_{\geq0}$ and set $s=m+3$. Define
\begin{equation}\label{eq:asymptoticVarianceConstant}
 \mathcal V_s=
 \frac{\Gamma(2s-1)\zeta(2s-2)}{\Gamma(s)^2}
 -\frac{s^2\zeta(s)^2}{2\zeta(2)}.
\end{equation}
Then $\mathcal V_s>0$ and, as $n\to\infty$,
\begin{align}
 \E[X_{m,n}]&\sim
 \zeta(s)\left(\frac{\sqrt{6n}}{\pi}\right)^s,
 \label{eq:asymptoticExpectation}\\
 \Var(X_{m,n})&\sim
 \mathcal V_s\left(\frac{\sqrt{6n}}{\pi}\right)^{2s-1}.
 \label{eq:asymptoticVariance}
\end{align}
In particular, the expectation has order $n^{(m+3)/2}$ and the
variance has order $n^{m+5/2}$.
\end{theorem}

The variance deserves special care. Since the expectation has order
$n^{s/2}$, its square and the second moment both have order $n^s$.
Their leading terms cancel in
$\Var(X_{m,n})=\E[X_{m,n}^2]-\E[X_{m,n}]^2$; the smaller
variance has order $n^{s-1/2}$. Consequently, leading-order
approximations to the two moments do not suffice. We must retain the
next term in the coefficient estimate, and it is this correction that
produces the second term in $\mathcal V_s$.

We isolate the required estimate in the following lemma. Its proof is
included to indicate precisely where the correction originates. As is
standard in saddle-point arguments, the uniform sectorial assumption
controls the factor being transferred near the saddle point, while
nonnegativity of its coefficients bounds it on the rest of the
integration circle.

\begin{lemma}\label{lemma:saddleTransfer}
Set $A=\zeta(2)$ and let $t_n>0$ be the unique solution of
\begin{equation}\label{eq:saddleDefinition}
 n=\sum_{k\geq1}\frac{k}{e^{kt_n}-1}.
\end{equation}
Suppose $F(x)$ is a power series with nonnegative coefficients and,
for some fixed $r>0$ and constants $c_0,c_1$, satisfies, uniformly
for $|\arg z|\leq\pi/4$ as $z\to0$,
\begin{equation}\label{eq:saddleHypothesis}
 F(e^{-z})=c_0z^{-r}+c_1z^{-r+1}
       +O\!\left(|z|^{-r+2}(1+|\log z|)\right).
\end{equation}
Then
\begin{equation}\label{eq:saddleTransfer}
 \frac{[x^n]P(x)F(x)}{p(n)}
 =c_0t_n^{-r}
 +\left(c_1+\frac{r(2-r)}{4A}c_0\right)t_n^{-r+1}
 +o(t_n^{-r+1}).
\end{equation}
If $c_0>0$ and only
$F(e^{-z})=c_0z^{-r}(1+o(1))$ is known uniformly in the same
sector, the corresponding leading-order conclusion is
$c_0t_n^{-r}(1+o(1))$.
\end{lemma}

\begin{proof}
Put $\Phi(z)=\log P(e^{-z})$, initially for $\Re z>0$. The
classical expansion of Euler's product and its differentiated forms
give, as $t\downarrow0$,
\begin{equation}\label{eq:logEulerExpansion}
 \Phi(t)=\frac{A}{t}+\frac12\log\frac{t}{2\pi}+O(t),
 \quad
 \Phi''(t)=\frac{2A}{t^3}(1+O(t)),\quad
 \Phi^{(3)}(t)=-\frac{6A}{t^4}(1+O(t)).
\end{equation}
For clarity, the first expansion is the logarithmic form of the
classical asymptotics of the Euler product; the derivative estimates
are part of that classical expansion and are not inferred merely by
differentiating an unspecified $O(t)$ term.

Differentiating the defining product before taking asymptotics gives
\[
 \Phi'(t)=-\sum_{k\geq1}\frac{k}{e^{kt}-1}.
\]
The sum on the right is continuous and strictly decreasing from
$+\infty$ to $0$ as $t$ increases from $0$ to $+\infty$. Thus
\eqref{eq:saddleDefinition} has a unique solution, characterized by
$\Phi'(t_n)=-n$. Moreover, since
$-\Phi'(t)\sim A/t^2$, we have
\begin{equation}\label{eq:tnAsymptotic}
 t_n\sim\sqrt{\frac{A}{n}}=\frac{\pi}{\sqrt{6n}}.
\end{equation}
This is the source of the scale $\sqrt{6n}/\pi$ in the theorem.

Write $t=t_n$ and $B=\Phi''(t)$. Cauchy's coefficient formula on
$|x|=e^{-t}$, with $x=e^{-t+i\theta}$ and
$dx=ix\,d\theta$, gives
\begin{equation}\label{eq:saddleIntegral}
 \frac{e^{nt+\Phi(t)}}{2\pi}
 \int_{-\pi}^{\pi}
 \exp\!\left(\Phi(t-i\theta)-\Phi(t)-in\theta\right)
 F(e^{-t+i\theta})\,d\theta
\end{equation}
for $[x^n]P(x)F(x)$. Replacing $F$ by $1$ gives the analogous
formula for $p(n)$. Their ratio therefore equals the ratio of the
two integrals: the common factor $e^{nt+\Phi(t)}/(2\pi)$ cancels.

We first explain why only a small interval around $\theta=0$
matters. The identity
\[
 \log\frac{|P(e^{-t+i\theta})|}{P(e^{-t})}
 =-\frac12\sum_{k\geq1}
 \log\!\left(1+\frac{2e^{-kt}(1-\cos k\theta)}
                      {(1-e^{-kt})^2}\right)
\]
follows from
$|1-ue^{i\alpha}|^2=(1-u)^2+2u(1-\cos\alpha)$.
For $|\theta|\leq t$, the summands with
$1\leq k\leq(2t)^{-1}$ show that the logarithm of the ratio is at
most $-c\theta^2/t^3$: in this range
$1-\cos(k\theta)\gg k^2\theta^2$ and
$1-e^{-kt}\asymp kt$, and there are $\asymp t^{-1}$ such terms.
For $t\leq|\theta|\leq\pi$, one instead uses
$t^{-1}\leq k\leq2t^{-1}$. In this range the factors depending on
$kt$ are bounded above and below, while the average of
$1-\cos(k\theta)$ over these consecutive indices is bounded away
from zero. The latter assertion follows, for $|\theta|\geq Ct$,
from the finite geometric-sum bound for $\sum e^{ik\theta}$; for
$t\leq|\theta|\leq Ct$, it follows by comparison with the
positive integral $\int_1^2(1-\cos(uv))\,du$, uniformly for
$1\leq|v|\leq C$. Consequently, for a fixed $c>0$,
\[
 \frac{|P(e^{-t+i\theta})|}{P(e^{-t})}
 \leq\begin{cases}
 \exp(-c\theta^2/t^3),&|\theta|\leq t,\\
 \exp(-c/t),&t\leq|\theta|\leq\pi.
 \end{cases}
\]
Since $F$ has nonnegative coefficients,
$|F(e^{-t+i\theta})|\leq F(e^{-t})=O(t^{-r})$.
The region outside
\[
 |\theta|\leq t^{3/2}\log(1/t)
\]
therefore contributes less than any fixed power of $t$ relative to
the Gaussian main terms. This interval is the \emph{major arc}; its
width is slightly larger than the natural Gaussian width $t^{3/2}$.
On the major arc $|\theta|/t\to0$, so the assumed sectorial
expansion for $F$ is available at $z=t-i\theta$.

Next, expand the exponent of \eqref{eq:saddleIntegral} about
$\theta=0$. Taylor's formula gives
\begin{align*}
 \Phi(t-i\theta)-\Phi(t)-in\theta
 &=-i\theta(\Phi'(t)+n)-\frac{\Phi''(t)\theta^2}{2}
   +\frac{i\Phi^{(3)}(t)\theta^3}{6}
   +O(t^{-5}\theta^4)\\
 &=-\frac{B\theta^2}{2}
   +\frac{i\Phi^{(3)}(t)\theta^3}{6}
   +O(t^{-5}\theta^4),
\end{align*}
where the linear term vanishes because $\Phi'(t)=-n$.
Since $B\sim2A/t^3$, the leading factor
$e^{-B\theta^2/2}$ is Gaussian, with typical
$\theta=O(t^{3/2})$. The terms involving $\theta^3$ and
$\theta^4$ in the exponent have respective sizes
$O(t^{1/2})$ and $O(t)$ on that scale.

The expansion \eqref{eq:saddleHypothesis} gives
\begin{align*}
 F(e^{-t+i\theta})
 &=c_0(t-i\theta)^{-r}+c_1(t-i\theta)^{-r+1}
       +O\!\left(t^{-r+2}(1+|\log t|)\right)\\
 &=c_0t^{-r}\left(1+\frac{ri\theta}{t}
      -\frac{r(r+1)}2\frac{\theta^2}{t^2}\right)
   +c_1t^{-r+1}+\text{higher-order terms}.
\end{align*}
Here we used
$(1-u)^{-r}=1+ru+\frac{r(r+1)}2u^2+O(u^3)$
with $u=i\theta/t$. A term odd in $\theta$ integrates to zero
against the leading Gaussian. The terms independent of $F$ in the
higher expansion of the exponential occur in both numerator and
denominator and cancel in their ratio. Thus, at relative order $t$,
the surviving corrections to the $c_0t^{-r}$ term come from the
quadratic term in $F$ and from the product of its linear term with
the cubic term of the exponent.

For the centered Gaussian density proportional to
$e^{-B\theta^2/2}$, the moments are
\[
 \E_{\mathrm G}[\theta^2]=B^{-1},\qquad
 \E_{\mathrm G}[\theta^4]=3B^{-2}.
\]
The quadratic term contributes
$-r(r+1)/(2Bt^2)$. The interaction of the two odd terms contributes
\[
 \frac{ri\theta}{t}\frac{i\Phi^{(3)}(t)\theta^3}{6}
 =-\frac{r\Phi^{(3)}(t)}{6t}\theta^4,
\]
whose Gaussian average is $-r\Phi^{(3)}(t)/(2tB^2)$.
Using \eqref{eq:logEulerExpansion}, their sum is
\begin{align*}
 -\frac{r(r+1)}{2Bt^2}
 -\frac{r\Phi^{(3)}(t)}{2tB^2}
 &=-\frac{r(r+1)}{4A}t+\frac{3r}{4A}t+o(t)\\
 &=\frac{r(2-r)}{4A}t+o(t).
\end{align*}
The remainder in \eqref{eq:saddleHypothesis} contributes only
$O(t^2(1+|\log t|))=o(t)$ relative to $c_0t^{-r}$ (when
$c_0\ne0$); the remaining Gaussian terms likewise contribute
$o(t)$ to the ratio. We conclude that
\[
 \frac{[x^n]P(x)F(x)}{p(n)}
 =c_0t^{-r}\left(1+\frac{r(2-r)}{4A}t\right)
   +c_1t^{-r+1}+o(t^{-r+1}),
\]
which is \eqref{eq:saddleTransfer}. Retaining only the leading
Gaussian term proves the final assertion.
\end{proof}

\begin{proof}[Proof of Theorem~\ref{theorem:asymptoticMoments}]
We first determine the singular expansions to which the saddle-point
lemma will be applied. Write $s=m+3$, so
\[
 a_k=\binom{k+s-2}{s-1}
 =\frac{k(k+1)\cdots(k+s-2)}{(s-1)!}.
\]
The leading coefficient of this polynomial in $k$ is
$1/(s-1)!=1/\Gamma(s)$. Its next coefficient is the sum of the
integers $0,1,\ldots,s-2$, divided by $(s-1)!$. Hence
\begin{equation}\label{eq:akAsymptotic}
 a_k=\frac{k^{s-1}}{\Gamma(s)}
      +\frac{s-2}{2\Gamma(s-1)}k^{s-2}
      +O(k^{s-3}).
\end{equation}
The lower-order terms are also polynomial in $k$; thus the Dirichlet
series associated with $a_k$ is a finite linear combination of
shifted zeta functions.

To expand $L_m$, substitute $x=e^{-t}$:
\[
 L_m(e^{-t})=\sum_{k\geq1}\frac{a_k}{e^{kt}-1}.
\]
For $\Re z>1$, the geometric series
$(e^u-1)^{-1}=\sum_{j\geq1}e^{-ju}$ gives
\[
 \int_0^\infty\frac{u^{z-1}}{e^u-1}\,du
 =\sum_{j\geq1}\int_0^\infty u^{z-1}e^{-ju}\,du
 =\Gamma(z)\zeta(z).
\]
Applying inverse Mellin transformation with $u=kt$ and interchanging
sum and integral on the absolutely convergent line $\Re z=c>s$
yields
\[
 L_m(e^{-t})=\frac{1}{2\pi i}
 \int_{c-i\infty}^{c+i\infty}
 \Gamma(z)\zeta(z)
 \left(\sum_{k\geq1}a_kk^{-z}\right)t^{-z}\,dz.
\]
By \eqref{eq:akAsymptotic},
\[
 \sum_{k\geq1}a_kk^{-z}
 =\frac{\zeta(z-s+1)}{\Gamma(s)}
 +\frac{s-2}{2\Gamma(s-1)}\zeta(z-s+2)
 +\text{lower shifts of }\zeta.
\]
The two rightmost poles occur at $z=s$ and $z=s-1$.
The simple pole of $\zeta(z-s+1)$ at $z=s$ has residue $1$;
therefore the residue of the full integrand at that point is
\[
 \Gamma(s)\zeta(s)\frac1{\Gamma(s)}t^{-s}
 =\zeta(s)t^{-s}.
\]
The residue at $z=s-1$ is similarly
\[
 \Gamma(s-1)\zeta(s-1)
 \frac{s-2}{2\Gamma(s-1)}t^{-(s-1)}
 =\frac{s-2}{2}\zeta(s-1)t^{-(s-1)}.
\]
Shifting the contour past these poles gives, uniformly for
$|\arg t|\leq\pi/4$,
\begin{equation}\label{eq:LMellin}
 L_m(e^{-t})=\zeta(s)t^{-s}
 +\frac{s-2}{2}\zeta(s-1)t^{-(s-1)}
 +O\!\left(|t|^{-(s-2)}(1+|\log t|)\right).
\end{equation}
The logarithm in the error allowance accommodates possible coincident
poles farther to the left. The contour shifts and their uniformity
follow from the exponential decay of $\Gamma(z)$ on vertical lines;
for $|\arg t|\leq\pi/4$, the growth of $t^{-z}$ does not overwhelm
this decay.

Next consider $Q_m(e^{-t})$. Write
\[
 w(u)=\frac{e^{-u}}{(1-e^{-u})^2}
      =\sum_{j\geq1}je^{-ju}.
\]
Its Mellin transform, for $\Re z>2$, is
\[
 \int_0^\infty u^{z-1}w(u)\,du
 =\Gamma(z)\sum_{j\geq1}j^{1-z}
 =\Gamma(z)\zeta(z-1).
\]
Consequently,
\[
 Q_m(e^{-t})=\frac1{2\pi i}
 \int_{c-i\infty}^{c+i\infty}
 \Gamma(z)\zeta(z-1)
 \left(\sum_{k\geq1}a_k^2k^{-z}\right)t^{-z}\,dz,
 \qquad c>2s-1.
\]
Since $a_k^2=k^{2s-2}/\Gamma(s)^2+O(k^{2s-3})$,
the rightmost pole of the Dirichlet factor is at $z=2s-1$, with
residue $1/\Gamma(s)^2$. The other factors are regular there, and
the resulting residue gives
\begin{equation}\label{eq:QMellin}
 Q_m(e^{-t})=
 \frac{\Gamma(2s-1)\zeta(2s-2)}{\Gamma(s)^2}
 t^{-(2s-1)}(1+o(1)).
\end{equation}
This estimate is uniform in the same sector.

Put
\[
 U=\zeta(s),\qquad
 V=\frac{s-2}{2}\zeta(s-1),\qquad
 C=\frac{\Gamma(2s-1)\zeta(2s-2)}{\Gamma(s)^2},
 \qquad A=\zeta(2).
\]
Equation \eqref{eq:LMellin} says
$L_m(e^{-t})=Ut^{-s}+Vt^{-s+1}+O(t^{-s+2}(1+|\log t|))$.
Apply Lemma~\ref{lemma:saddleTransfer} with $F=L_m$, $r=s$,
$c_0=U$, and $c_1=V$. By \eqref{eq:expectation},
\begin{equation}\label{eq:meanTwoTerm}
 \E[X_{m,n}]
 =Ut_n^{-s}
 +\left(V+\frac{s(2-s)}{4A}U\right)t_n^{-s+1}
 +o(t_n^{-s+1}).
\end{equation}
In particular, $\E[X_{m,n}]\sim Ut_n^{-s}$, which already proves
the claimed expectation asymptotic after using
\eqref{eq:tnAsymptotic}.

For the variance, first square the expansion of $L_m$:
\[
 L_m(e^{-t})^2
 =U^2t^{-2s}+2UVt^{-2s+1}
 +O\!\left(t^{-2s+2}(1+|\log t|)\right).
\]
Applying the saddle-point lemma to $F=L_m^2$ with $r=2s$
therefore gives
\[
 \frac{[x^n]P(x)L_m(x)^2}{p(n)}
 =U^2t_n^{-2s}
 +\left(2UV+\frac{(2s)(2-2s)}{4A}U^2\right)
 t_n^{-(2s-1)}+o(t_n^{-(2s-1)}).
\]
The leading-order version of the lemma applied to
\eqref{eq:QMellin} supplies the remaining term,
\[
 \frac{[x^n]P(x)Q_m(x)}{p(n)}
 =Ct_n^{-(2s-1)}(1+o(1)).
\]
By \eqref{eq:secondMoment}, the sum of these two expressions is
\begin{equation}\label{eq:secondMomentTwoTerm}
 \E[X_{m,n}^2]
 =U^2t_n^{-2s}
 +\left(2UV+\frac{(2s)(2-2s)}{4A}U^2+C\right)
   t_n^{-(2s-1)}
 +o(t_n^{-(2s-1)}).
\end{equation}

We now subtract the square of \eqref{eq:meanTwoTerm} from
\eqref{eq:secondMomentTwoTerm}. The square of the mean is
\[
 \E[X_{m,n}]^2
 =U^2t_n^{-2s}
 +\left(2UV+\frac{2s(2-s)}{4A}U^2\right)
 t_n^{-(2s-1)}
 +o(t_n^{-(2s-1)}).
\]
The $U^2t_n^{-2s}$ terms cancel, as do the $2UV$ terms. The
remaining saddle corrections have coefficient
\[
 \frac{(2s)(2-2s)}{4A}
 -\frac{2s(2-s)}{4A}
 =-\frac{s^2}{2A}.
\]
We therefore obtain
\begin{align*}
 \Var(X_{m,n})
 &=\left(C-\frac{s^2U^2}{2A}\right)t_n^{-(2s-1)}
 +o(t_n^{-(2s-1)})\\
 &=\mathcal V_s\,t_n^{-(2s-1)}
 +o(t_n^{-(2s-1)}).
\end{align*}
This calculation explains why the contribution from $Q_m$ alone is
not the correct variance constant: fixing the size $n$ imposes a
constraint that reduces the fluctuations. Once positivity of
$\mathcal V_s$ is established, \eqref{eq:tnAsymptotic} gives
\eqref{eq:asymptoticVariance}.

Finally, we prove $\mathcal V_s>0$. With
$w(u)=e^{-u}(1-e^{-u})^{-2}$, the expansion
$w(u)=\sum_{j\geq1}je^{-ju}$ gives, for $a>1$,
\[
 \int_0^\infty u^a w(u)\,du
 =\Gamma(a+1)\zeta(a).
\]
Define $f(u)=u^{s-1}/\Gamma(s)$ and $h(u)=u$. Then
\begin{align*}
 \int_0^\infty f(u)^2w(u)\,du&=C,\\
 \int_0^\infty f(u)h(u)w(u)\,du&=sU,\\
 \int_0^\infty h(u)^2w(u)\,du&=2A.
\end{align*}
The Cauchy--Schwarz inequality in the weighted space
$L^2((0,\infty),w(u)\,du)$ gives $(sU)^2\leq(2A)C$, or
$C-s^2U^2/(2A)\geq0$. Equality would require
$u^{s-1}/\Gamma(s)$ to be proportional to $u$ almost everywhere.
That is impossible for $s\geq3$, so the inequality is strict:
$\mathcal V_s>0$. This completes the proof.
\end{proof}

\begin{remark}\label{rem:varianceCorrection}
The negative contribution $-s^2\zeta(s)^2/(2\zeta(2))$ has a
probabilistic interpretation. In the independent-geometric model of
random partitions, let $N_k$ be independent geometric random variables
with law $\Pr(N_k=j)=(1-e^{-kt})e^{-jkt}$ for $j\geq0$,
and put
$T=\sum_{k\geq1}kN_k$ and $Y=\sum_{k\geq1}a_kN_k$.
Conditioning on $T=n$ gives the uniform distribution on partitions of
$n$. In the unconditioned model,
$\Var(Y)=Q_m(e^{-t})$,
$\Var(T)\sim2At^{-3}$, and
$\operatorname{Cov}(Y,T)\sim sUt^{-s-1}$.
The Gaussian heuristic for conditioning on $T=n$ subtracts
$\operatorname{Cov}(Y,T)^2/\Var(T)$ from the variance of $Y$;
this has leading constant $s^2U^2/(2A)$. This heuristic agrees
with the rigorous saddle-point calculation, but the preceding proof
does not rely on a conditional central limit theorem.
\end{remark}

We now apply the theorem to the normalized statistics introduced in
Section~\ref{sec:moments}. These applications illustrate two different asymptotic regimes:
under uniform sampling the normalized Gini statistic tends to its
upper endpoint, whereas the min--max normalized Wiener index tends
to its lower endpoint, despite having the same extremal partitions.

\begin{corollary}\label{cor:normalizedGiniAsymptotic}
Let $\widetilde g(\lambda)=g(\lambda)/\binom n2$ for
$\lambda\vdash n\geq2$. Put $\mathcal V_3$ as in
\eqref{eq:asymptoticVarianceConstant}. Then
\begin{align}
 1-\E[\widetilde g(\Lambda_n)]
 &\sim 2\zeta(3)\left(\frac{\sqrt6}{\pi}\right)^3n^{-1/2},
 \label{eq:giniMeanAsymptotic}\\
 \Var(\widetilde g(\Lambda_n))
 &\sim 4\mathcal V_3\left(\frac{\sqrt6}{\pi}\right)^5n^{-3/2}.
 \label{eq:giniVarianceAsymptotic}
\end{align}
In particular, $\widetilde g(\Lambda_n)\to1$ in probability.
\end{corollary}
\begin{proof}
From \eqref{gini} and \eqref{eq:defineMm},
\[
 \widetilde g(\lambda)
 =\frac{n+1}{n-1}-\frac{M_0(\lambda)}{\binom n2}.
\]
Hence
\[
 1-\E[\widetilde g(\Lambda_n)]
 =-\frac2{n-1}+\frac{\E[X_{0,n}]}{\binom n2},
 \qquad
 \Var(\widetilde g(\Lambda_n))
 =\frac{\Var(X_{0,n})}{\binom n2^2}.
\]
Set $m=0$, and thus $s=3$, in
Theorem~\ref{theorem:asymptoticMoments}. Since
$\binom n2\sim n^2/2$, the expectation contribution has leading
order $n^{3/2}/n^2=n^{-1/2}$, whereas $2/(n-1)=O(n^{-1})$
is smaller. The variance contribution has order
$n^{5/2}/n^4=n^{-3/2}$. These comparisons give the two displayed
asymptotics. Finally, because $0\leq\widetilde g\leq1$,
Markov's inequality applied to the nonnegative variable
$1-\widetilde g(\Lambda_n)$ yields, for every $\varepsilon>0$,
\[
 \Pr\!\left(1-\widetilde g(\Lambda_n)>\varepsilon\right)
 \leq\frac{1-\E[\widetilde g(\Lambda_n)]}{\varepsilon}
 \longrightarrow0.
\]
Thus $\widetilde g(\Lambda_n)\to1$ in probability.
\end{proof}

\begin{corollary}\label{cor:normalizedWienerAsymptotic}
For $n\geq3$, let $\widetilde W$ be the min--max normalized Wiener
index defined in Section~\ref{sec:moments}. Then
\begin{equation}\label{eq:wienerMeanAsymptotic}
 \E[\widetilde W(\Lambda_n)]
 \sim6\zeta(3)\left(\frac{\sqrt6}{\pi}\right)^3n^{-1/2}.
\end{equation}
In particular, $\widetilde W(\Lambda_n)\to0$ in probability.
\end{corollary}
\begin{proof}
Proposition~\ref{prop:wienerWeightedRelation}, together with the
size-preserving bijection $\lambda\mapsto\lambda^*$, gives
\[
 \E[W(S_{\Lambda_n})]
 =(n+2)\E[X_{0,n}]-2\E[X_{1,n}].
\]
By \eqref{eq:asymptoticExpectation},
$\E[X_{0,n}]\sim\zeta(3)(\sqrt6/\pi)^3n^{3/2}$,
whereas $\E[X_{1,n}]=O(n^2)$. The first term therefore has
order $n^{5/2}$ and dominates the second, so
\[
 \E[W(S_{\Lambda_n})]
 \sim\zeta(3)\left(\frac{\sqrt6}{\pi}\right)^3n^{5/2}.
\]
Subtracting $n^2$ does not alter this leading term. Since
$\binom n3\sim n^3/6$, the definition
$\widetilde W=(W-n^2)/\binom n3$ yields
\eqref{eq:wienerMeanAsymptotic}. Finally,
$0\leq\widetilde W\leq1$ and its expectation tends to $0$;
Markov's inequality gives convergence to $0$ in probability.
\end{proof}

For $m=0$ the variance constant has the more explicit form
\[
 \mathcal V_3
 =6\zeta(4)-\frac{9\zeta(3)^2}{2\zeta(2)}
 =\frac{\pi^4}{15}-\frac{27\zeta(3)^2}{\pi^2}>0.
\]
The asymptotic results describe not only the magnitude of a typical
weighted sum but also its concentration: the standard deviation of
$X_{m,n}$ has order $n^{(2s-1)/4}$, while its mean has order
$n^{s/2}$, so their ratio decays at order $n^{-1/4}$ for fixed $m$.

\section{Dominance and Antichains}\label{sec:dominanceAntichains}

A second application of the weighted sums concerns the partial order
on partitions defined by comparing initial row sums. For
$\lambda,\mu\vdash n$, write $\lambda\preceq\mu$ when
\begin{equation}\label{eq:dominanceOrder}
 \sum_{i=1}^r\lambda_i\leq\sum_{i=1}^r\mu_i
 \qquad(1\leq r\leq n),
\end{equation}
where partitions are padded by zeros. Write $\lambda\prec\mu$ if
$\lambda\preceq\mu$ and $\lambda\ne\mu$. This is the
\emph{dominance order}; a partition higher in this order places more
of its total size in its early rows. For example,
$(1^n)\preceq\lambda\preceq(n)$ for every $\lambda\vdash n$.
The partitions of $n$ form a lattice $\mathcal L_n$ under dominance
\cite{Brylawski:1,Inequalities:1}.

An \emph{antichain} in a partially ordered set is a subset in which
no two distinct elements are comparable. The \emph{width} $w(n)$ of
$\mathcal L_n$ is the maximum possible cardinality of such a subset;
see \cite{Antichain:1,OEIS:A076269}. Finding the width can be
difficult, so explicit families of antichains give useful lower
bounds. The central observation here is that $M_m$ is strictly
\emph{decreasing} with dominance for every $m\geq0$. Therefore,
two different partitions having the same value of $M_m$ cannot be
comparable.

\begin{proposition}\label{prop:dominancelattice}
Let $m\in\mathbb Z_{\geq0}$ and let $\lambda\prec\mu$ be distinct
partitions of $n$. Then
\begin{equation}\label{eq:dominanceWeighted}
 M_m(\lambda)>M_m(\mu)
\end{equation}
and
\begin{equation}\label{eq:dominanceBinomial}
 \sum_{i=1}^n\binom{\lambda_i+m+1}{m+2}
 <\sum_{i=1}^n\binom{\mu_i+m+1}{m+2}.
\end{equation}
\end{proposition}
\begin{proof}
Let $w_i=\binom{i+m}{m+1}$ denote the weight in row $i$.
Define the differences of initial row sums by
\[
 D_r=\sum_{i=1}^r(\mu_i-\lambda_i),\qquad
 D_0=0.
\]
Dominance gives $D_r\geq0$ for each $r$. Moreover,
$D_n=0$, since both partitions have total size $n$. Because the
partitions are distinct, at least one of $D_1,\ldots,D_{n-1}$
is strictly positive. Otherwise all successive differences
$D_r-D_{r-1}=\mu_r-\lambda_r$ would vanish, forcing
$\lambda=\mu$.

We now use a discrete summation-by-parts identity. As
$\mu_i-\lambda_i=D_i-D_{i-1}$, we have
\begin{align*}
 M_m(\lambda)-M_m(\mu)
 &=-\sum_{i=1}^nw_i(D_i-D_{i-1})\\
 &=-\sum_{i=1}^nw_iD_i+\sum_{i=1}^nw_iD_{i-1}\\
 &=\sum_{r=1}^{n-1}(w_{r+1}-w_r)D_r.
\end{align*}
The boundary terms vanish because $D_0=D_n=0$. By Pascal's
identity,
\[
 w_{r+1}-w_r
 =\binom{r+m+1}{m+1}-\binom{r+m}{m+1}
 =\binom{r+m}{m}>0.
\]
Consequently,
\[
 M_m(\lambda)-M_m(\mu)
 =\sum_{r=1}^{n-1}\binom{r+m}{m}D_r>0,
\]
proving \eqref{eq:dominanceWeighted}. The strictness comes from
having at least one positive $D_r$ and \emph{every} weight increment
positive; merely nondecreasing weights would not ensure strictness.

Conjugation reverses strict dominance
\cite{Brylawski:1}: $\lambda\prec\mu$ implies
$\mu^*\prec\lambda^*$. Applying the first inequality to this
reversed pair gives
$M_m(\mu^*)>M_m(\lambda^*)$. Finally,
Proposition~\ref{prop:generalcase} shows that
\[
 M_m(\nu^*)=\sum_i\binom{\nu_i+m+1}{m+2}.
\]
Applying this with $\nu=\lambda$ and $\nu=\mu$ proves
\eqref{eq:dominanceBinomial}.
\end{proof}

The two inequalities offer complementary viewpoints. The row-index
sum $M_m(\lambda)$ decreases with dominance, since shifting weight
into earlier rows tends to replace larger row indices by smaller
ones. Conversely, the sum of binomial functions of the part sizes
increases with dominance. The latter is consistent with the strict
convexity of $k\mapsto\binom{k+m+1}{m+2}$ on the nonnegative
integers and with the usual majorization inequalities.

\begin{corollary}\label{cor:dominanceAntichains}
For $m,r\in\mathbb Z_{\geq0}$, the level set
\begin{equation}\label{eq:levelSet}
 \mathcal A_{m,r}(n)=\{\lambda\vdash n:M_m(\lambda)=r\}
\end{equation}
is an antichain in $\mathcal L_n$.
\end{corollary}
\begin{proof}
Suppose distinct members $\lambda,\mu$ of this level set were
comparable. One would strictly dominate the other, and
Proposition~\ref{prop:dominancelattice} would imply that their
$M_m$ values were strictly different. This contradicts the
defining equality of the level set.
\end{proof}

Because the bivariate generating function counts the partitions in
each level set, the antichain construction translates immediately into
coefficient bounds. The parameter $m$ provides a family of different
ways to slice the same dominance lattice into antichains.

\begin{corollary}\label{cor:dominanceWidth}
For every $m\in\mathbb Z_{\geq0}$ and $n\geq1$,
\begin{equation}\label{eq:dominanceWidthBound}
 w(n)\geq\max_{r\geq0}|\mathcal A_{m,r}(n)|
 =\max_{r\geq0}[q^rx^n]H_m(q,x).
\end{equation}
\end{corollary}
\begin{proof}
Every $\mathcal A_{m,r}(n)$ is an antichain by
Corollary~\ref{cor:dominanceAntichains}, so its cardinality cannot
exceed the maximum antichain size $w(n)$. By the definition of
$H_m(q,x)$, the coefficient $[q^rx^n]H_m(q,x)$ counts precisely
those partitions $\lambda\vdash n$ with $M_m(\lambda)=r$.
Taking the maximum over all possible $r$ gives the claim.
\end{proof}

The $m=0$ construction is related to the Gini-index level sets in
\cite{Kopitzke:1}, because
$g(\lambda)=\binom{n+1}{2}-M_0(\lambda)$. Thus equality of
$M_0$ values is equivalent to equality of Gini values for a fixed
$n$. The corresponding lower-bound sequence appears as
\seqnum{A337206} in the \emph{On-Line Encyclopedia of Integer
Sequences} \cite{OEIS:1}; the widths themselves appear as
\seqnum{A076269}. Although each $m$ supplies a valid lower bound,
there is no reason for those bounds to increase as $m$ increases:
different choices of weights can separate or identify different
incomparable partitions.

\begin{example}\label{ex:incomparable}
Consider $\lambda=(6,4)$ and $\mu=(8,1,1)$, two partitions of
$10$. They are incomparable because the first initial row sums
satisfy $6<8$, whereas the first two satisfy $6+4>8+1$.
The first four statistics give
\[
 \begin{array}{c|rrrr}
 m&0&1&2&3\\ \hline
 M_m(\lambda)&14&18&22&26\\
 M_m(\mu)&13&17&22&28
 \end{array}
\]
The $m=2$ values are equal, so the pair is certified as
incomparable by the preceding proposition: distinct comparable
partitions cannot have the same $M_2$ value. Alternatively,
$M_0(\lambda)>M_0(\mu)$ would force the direction
$\lambda\prec\mu$ if the pair were comparable, while
$M_3(\lambda)<M_3(\mu)$ would force the opposite direction.
The conflicting inequalities therefore rule out comparability.
An inequality for a \emph{single} $M_m$ is not enough to infer a
dominance relation, since the converse of
Proposition~\ref{prop:dominancelattice} need not hold.
\end{example}

As a more concrete example of a full level set, the three partitions
\[
 (7,1,1,1),\qquad(6,2,2),\qquad(5,4,1)
\]
of $10$ all have $M_0=16$. They therefore form an antichain,
without requiring pairwise tests of the dominance inequalities.
Other level sets of $M_0$ on partitions of $10$ also have three
members.

\begin{table}[htbp]
\centering
\begin{tabular}{c@{\qquad}ccccc}
\toprule
$m$ & $0$ & $1$ & $2$ & $3$ & $4$ \\
\midrule
$\max_r|\mathcal A_{m,r}(10)|$ & $3$ & $2$ & $2$ & $1$ & $2$ \\
\bottomrule
\end{tabular}
\caption{Level-set lower bounds for $w(10)$, obtained by enumerating
partitions of $10$. The actual width is $w(10)=4$
\cite{OEIS:A076269}.}
\label{tab:width10}
\end{table}

Table~\ref{tab:width10} illustrates two limitations of this method:
its bounds need not be monotonic in $m$, and none of those listed
achieves the true width for $n=10$. Nevertheless, the bounds are
computable directly from the coefficients of $H_m$, and varying $m$
can reveal incomparability that is not detected by another member
of the family. Understanding which weight parameters yield the
strongest bounds, and how close these bounds can come to $w(n)$,
remains a natural direction for further work.

\section{Concluding Remarks}\label{sec:conclusion}

The binomially weighted statistics $M_m$ admit several complementary
descriptions. A Young diagram makes them row-weighted cell sums;
conjugation converts them into sums of polynomial functions of part
sizes; and Frobenius coordinates express them as contributions from
principal hooks. The fact that these descriptions agree is the
combinatorial foundation of the generating-function results.

The bivariate generating function retains the distribution of each
statistic, not merely its sum over partitions of a fixed integer.
Differentiating that function gives exact formulas for the moments,
while its singular behavior and a refined saddle-point transfer
estimate lead to asymptotics for the expectation and variance. The
variance calculation is especially sensitive to the constraint that
the sampled partition have total size exactly $n$: the leading
second-moment contributions cancel, leaving a smaller variance with
a nontrivial positive constant.

The motivating Gini and Wiener examples demonstrate the reach of
the common weighted-sum framework: the former is determined by
$M_0(\lambda)$, and the latter by a linear combination of
$M_0(\lambda^*)$ and $M_1(\lambda^*)$. In particular, the same
moment asymptotics yield opposite limiting behavior under the
natural normalizations considered here: the normalized Gini statistic
converges to $1$ in probability, while the min--max normalized Wiener
index converges to $0$. The necessary conjugation in the Wiener
identity also illustrates why formulas for individual partitions
must be distinguished from identities valid only after averaging.

Finally, strict monotonicity under dominance turns every level set
of $M_m$ into an antichain, yielding an explicit family of lower
bounds for the width of the dominance lattice. Further questions
include limiting distributions of appropriately centered and scaled
weighted sums, higher-order moment asymptotics, and quantitative
comparisons of these lower bounds for different $m$. The bounds are
valid for every $m$, but when they approximate the true width well
is not yet understood.

\bigskip
\noindent\textit{2020 Mathematics Subject Classification:}
Primary 05A17; Secondary 05A15, 05E05, 06A06, 91B99.\\
\noindent\textit{Keywords:} integer partition, weighted sum,
generating function, Frobenius notation, dominance order, antichain.\\
\noindent\textit{OEIS sequences:} \seqnum{A076269}, \seqnum{A337206}.


\begin{thebibliography}{99}

\bibitem{Andrews:1}
George~E. Andrews and Kimmo Eriksson.
\newblock {\em Integer Partitions}.
\newblock Cambridge University Press, Cambridge, 2004.

\bibitem{Brylawski:1}
Thomas Brylawski.
\newblock The lattice of integer partitions.
\newblock {\em Discrete Mathematics}, 6(3):201--219, 1973.

\bibitem{Origins:1}
L.~Ceriani and P.~Verme.
\newblock The origins of the {G}ini index: Extracts from variabilit\`a e
  mutabilit\`a (1912) by {C}orrado {G}ini.
\newblock {\em Journal of Economic Inequality}, 10:1--23, 2012.

\bibitem{Antichain:1}
Edward Early.
\newblock Chain lengths in the dominance lattice.
\newblock {\em Discrete Mathematics}, 313(20):2168--2177, 2013.

\bibitem{Farris:1}
Frank~A. Farris.
\newblock The {G}ini index and measures of inequality.
\newblock {\em American Mathematical Monthly}, 117(10):851--864, 2010.

\bibitem{FlajoletSedgewick:1}
Philippe Flajolet and Robert Sedgewick.
\newblock {\em Analytic Combinatorics}.
\newblock Cambridge University Press, 2009.

\bibitem{Fristedt:1}
Bert Fristedt.
\newblock The structure of random partitions of large integers.
\newblock {\em Transactions of the American Mathematical Society},
  337(2):703--735, 1993.

\bibitem{Fulton:1}
William Fulton.
\newblock {\em Young Tableaux}.
\newblock London Mathematical Society Student Texts. Cambridge University
  Press, 1997.

\bibitem{WienerIndex:1}
Ivan Gutman, Juan Rada, and O.~Araujo.
\newblock The wiener index of starlike trees and a related partial order.
\newblock {\em MATCH Communications in Mathematical and in Computer Chemistry},
  42:145--154, 2000.

\bibitem{Kopitzke:1}
Grant Kopitzke.
\newblock The gini index of an integer partition.
\newblock {\em Journal of Integer Sequences}, 23, 2020.
\newblock Article 20.9.7.

\bibitem{Kopitzke:Dissertation}
Grant Kopitzke.
\newblock {\em The Gini Index in Algebraic Combinatorics and Representation
  Theory}.
\newblock PhD thesis, University of Wisconsin--Milwaukee, 2021.
\newblock Available at \url{https://arxiv.org/abs/2201.12867}.

\bibitem{Lorenz:1}
M.~O. Lorenz.
\newblock Methods of measuring the concentration of wealth.
\newblock {\em Publications of the American Statistical Association},
  9(70):209--219, 1905.

\bibitem{Inequalities:1}
Albert~W. Marshall, Ingram Olkin, and Barry~C. Arnold.
\newblock {\em Inequalities: Theory of Majorization and Its Applications}.
\newblock Springer Series in Statistics. Springer, 2 edition, 2011.

\bibitem{OEIS:A076269}
{OEIS Foundation Inc.}
\newblock {A076269}: Size of largest antichain in partition lattice {Par(n)}.
\newblock \url{https://oeis.org/A076269}, 2026.
\newblock Accessed November 11, 2026.

\bibitem{OEIS:1}
{OEIS Foundation Inc.}
\newblock The on-line encyclopedia of integer sequences.
\newblock \url{https://oeis.org}, 2026.
\newblock Accessed 2026-11-11.

\bibitem{Wilf:1}
Herbert~S. Wilf.
\newblock {\em Generatingfunctionology}.
\newblock Academic Press, 2 edition, 1994.

\end{thebibliography}
\end{document}